\documentclass[a4paper,11pt,oneside]{article}
\usepackage[english]{babel}

\usepackage[dvips]{graphicx}
\usepackage[svgnames]{xcolor}
\colorlet{foldline}{Red}
\colorlet{crease}{Gray}
\colorlet{halfgray}{DarkGrey}

\usepackage{amssymb}
\usepackage[tbtags]{amsmath}
\usepackage{amsthm}
\usepackage{icomma}

\theoremstyle{definition}
\newtheorem{definition}{Definition} 

\theoremstyle{plain}

\usepackage[T1]{fontenc}
\usepackage[osf,sc]{mathpazo}
\usepackage[euler-digits,euler-hat-accent]{eulervm}

\usepackage[top=32mm]{geometry}

\usepackage{caption}
\usepackage{booktabs}
\usepackage[flushleft]{threeparttable}

\usepackage{tabularx}

\usepackage[authoryear,round,longnamesfirst]{natbib}
\usepackage{fancyhdr}
\usepackage{hyperref}
\hypersetup{
  breaklinks = true,
  linktocpage=true,
  linkcolor  = MidnightBlue,
  citecolor  = MidnightBlue,
  urlcolor   = MidnightBlue,
  colorlinks = true,
}
\usepackage{breakcites}    

\usepackage{pst-all} 
\newpsobject{showgrid}{psgrid}{subgriddiv=1,griddots=10,gridlabels=6pt}
\psset{arrowscale=2}

\title{Geometric solution of a quintic equation by two-fold origami}
\author{Jorge C. Lucero\thanks{Dept.\ Computer Science, University of Bras\'{i}lia, Brazil. E-mail: \href{mailto:lucero@unb.br}{lucero@unb.br} }} 
\date{\today}

\begin{document}
\maketitle
\begin{abstract} 
\noindent
This article shows how to find the solution of an arbitrary quintic equation by performing two simultaneous folds on a sheet of paper. The folds achieve specific incidences between a set of points and lines that are determined by the coefficients of the quintic. Complete equations for computing the set are given, and their application is illustrated with an example.
\end{abstract}

\section{Introduction} 

It is well-known that arbitrary cubic equations may be geometrically solved by using origami (paper folding) \citep{Geretschlager1995, Hull2011}. In this technique, a sheet of paper is folded along a straight line so as to achieve a set of specific incidences between points and lines \citep{Lucero2017}, and these points and lines are determined by the coefficients of the equation. As a consequence, several classical construction problems related to cubic equations may be solved by sequences of folds, such as trisecting an arbitrary angle \citep{Hull1996}, duplicating the cube \citep{Messer1986}, and constructing the regular heptagon \citep{Geretschlaeger1997}. Also arbitrary quartic equations may be solved in this way after reduction to lower degree equations, although such a process may be excessively lengthy and complex for practical use \citep{Edwards2001}.  

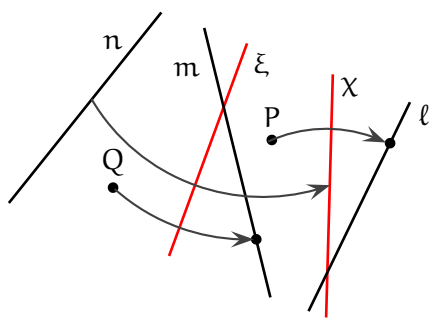
\begin{figure}[!htb]
\centering
\begin{pspicture}(-3,-2.5)(3.,2.5)
%\showgrid
\psset{xunit=.6cm,yunit=.6cm}
\rput{-20}{%
\psplot[linewidth=1pt]{-3.7}{-2}{x 3 add 3 mul}
\uput[180]{*0}(-2.3,2.2){$n$}
\psline[linewidth=1pt,linecolor=foldline](0,-2)(0,3)
\uput[0]{*0}(0,2.5){$\xi$}
\psplot[linewidth=1pt,linecolor=foldline]{2}{3.7}{x -3 add -3 mul}
\uput[0]{*0}(2,2.7){$\chi$}
\psplot[linewidth=1pt]{-1}{2.4}{x 1 add -1.5 mul 3 add}
\uput[180]{*0}(-.5,2.2){$m$}
\qdisk(-1.67,-1){2pt}
\uput[90]{*0}(-1.67,-1){$Q$}
\qdisk(1.67,-1){2pt}
\psplot[linewidth=1pt]{3.29}{3.8}{x -3.5 add 10 mul}
\uput[0]{*0}(3.8,2.7){$\ell$}
\qdisk(3.7,2){2pt}
\qdisk(1.24,1.18){2pt}
\uput[90]{*0}(1.24,1.18){$P$}
\psarcn[linecolor=darkgray]{->}(3.3,-1){1.8}{135}{84}
\psarc[linecolor=darkgray]{->}(0,3){2.6}{248}{292}
\psarc[linecolor=darkgray]{->}(0,4){2.6}{230}{310}
}
\end{pspicture}
\caption{Two-fold operation AL4a6ab. Lines $\chi$ and $\xi$ are the fold lines.} 
\label{nishi}
\end{figure}

In the above constructions, a single fold is performed at each step. If simultaneous folds are allowed, then solutions to higher degree equations may be constructed. A theorem by \citet{Alperin2006} states that every polynomial equation of degree $n$ with real solutions may be solved by performing $n-2$ simultaneous folds. In case of a quintic equation, the theorem demands three simultaneous folds. However, a recent paper by \citet{Nishimura2015} showed that the following two-fold operation is enough (the name ``AL4a6ab'' follows \citeauthor{Alperin2006}'s notation):

\begin{definition}[Operation AL4a6ab] Given two points $P$, $Q$ and three lines $\ell$, $m$, $n$, simultaneously fold along a line $\chi$ to place $P$ onto $\ell$, and a line $\xi$ to place $Q$ onto $m$ and to align $n$ and $\chi$ (see Fig.~\ref{nishi}).
\end{definition}

In his analysis, \citeauthor{Nishimura2015} assumed a coordinate plane such that $Q(0,1)$ and $m: y=-1$, and defined $P(p, q)$, $\ell: x= -k$ and $n: ax+by+1=0$, where $a$, $b$, $k$, $p$ and $q$ are real numbers. Then, he showed that the inclination of the fold line $\xi$ is a solution of  $t^5+\alpha t^4 + \beta t^3+\gamma t^2 +\delta t+\epsilon=0$, where coefficients $\alpha$, $\beta$, $\gamma$, $\delta$ and $\epsilon$ are also real. Next, he proved that appropriate values of $a$, $b$, $k$, $p$ and $q$ may be computed from the coefficients of a given quintic equation in depressed form (i.e., with $\alpha=0$), and under the condition $\epsilon^2-4(\beta + \delta +1)\ge 0$. In case of a general quintic, its depressed form is obtained by application of the Cardano transformation $t=t'-\frac{\alpha}{5}$, and the later condition may be met by a suitable scale change $t=ct'$.
  
This is a relevant result for the field of origami mathematics because it may simplify geometric constructions related to quintic equations, such as the quintisection of an arbitrary angle \citep{Lang2004a} and the construction of the regular hendecagon \citep{Lucerotoappear}. In order to perform two simultaneous folds, the folder must smoothly move both folds changing their position and orientation until all the required incidence constraints are met. With more than two simultaneous folds, such a manipulation becomes prohibitively difficult. However, practical applications of \citeauthor{Nishimura2015}'s \citeyearpar{Nishimura2015} analysis may find limitations because of the initial transformation into depressed form, and also because the analysis demands a line $n$ that does not to pass through the coordinate origin $(0, 0)$ (i.e., the analysis might miss a solution with a line $n$ through the origin, as it will be shown later with an example). In addition, the purpose of the analysis was to prove the existence of a geometrical solution, but the solution was not cast in explicit form. 

In this article, the above limitations will be solved in order to widen the applicability of the folding operation, and full equations for $P$, $Q$, $\ell$, $m$ and $n$ in terms of the given quintic coefficients will be provided.

\section{Analysis of operation AL4a6ab}

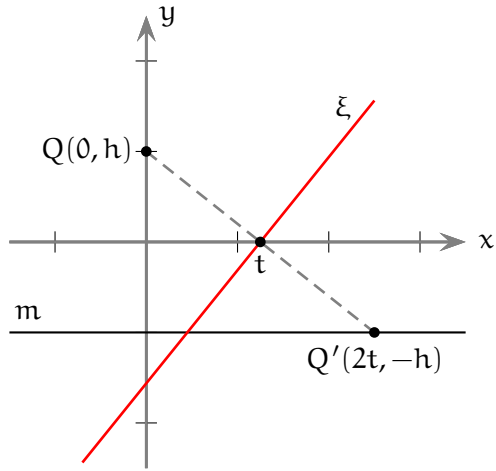
\begin{figure}[!htb]
\centering
\begin{pspicture}(-3,-3.5)(4,3)
\psset{xunit=1.2cm,yunit=1.2cm}
%\showgrid
\psaxes[linecolor=gray,labels=none,linewidth=1pt]{->}(0,0)(-1.5,-2.5)(3.5,2.5)[$x$,0][$y$,0]
\uput[180]{*0}(0,1){$Q(0, h)$}
\psline(-1.5,-1)(3.5,-1)
\uput[90]{*0}(-1.3,-1){$m$}
\uput[-90]{*0}(2.5,-1){$Q'(2t, -h)$}
\psline[linecolor=gray,linewidth=1pt, linestyle=dashed](0,1)(2.5,-1)
\psplot[linecolor=foldline,linewidth=1pt,plotpoints=10]{-.7}{2.5}{x -1.25 add 1.25 mul}
\uput[-90]{*0}(1.25,0){$t$}
\uput[90]{*0}(2.15,1.2){$\xi$}
\qdisk(0,1){2pt}
\qdisk(2.5,-1){2pt}
%\qdisk(2.5,1.56){2pt}
\qdisk(1.25,0){2pt}
\end{pspicture}
\caption{Fold line $\xi$ places point $Q$ onto line $m$ at $Q'$.} 
\label{PaLa}
\end{figure}

Assume a coordinate plane such that $Q(0, h)$ and $m: y=-h$ (Fig.~\ref{PaLa}).  

Fold $\xi$ places $Q$ onto line $m$ at $Q'(2t, -h)$, where $t$ is a parameter. The fold line may be described by a vector equation of the form $\mathbf{n}\cdot\mathbf{x}=\mathbf{n}\cdot\mathbf{x}_0$, where $\mathbf{n}$ is a normal vector, $\mathbf{x}=(x, y)$ and $\mathbf{x}_0$ is a point on the line. A normal vector to $\xi$ is $\overrightarrow{QQ'}=(2t, -2h)$ and $\xi$ passes through the midpoint of $\overline{QQ'}$ at $(t, 0)$. Therefore, $\xi$ has an equation
\begin{equation}
tx-hy = t^2.
\label{Fa}
\end{equation}

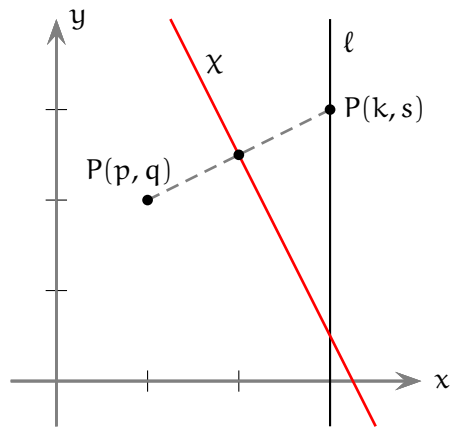
\begin{figure}[!htb]
\centering
\begin{pspicture}(-1,-1)(5,5.5)
\psset{xunit=1.2cm,yunit=1.2cm}
%\showgrid
\psaxes[linecolor=gray,labels=none,linewidth=1pt]{->}(0,0)(-.5,-.5)(4,4)[$x$,0][$y$,0]
\uput[90]{0}(.8,2){$P(p, q)$}
\psline(3,-.5)(3,4)
\uput[0]{*0}(3,3.75){$\ell$}
\uput[0]{*0}(3,3){$P(k, s)$}
\psline[linecolor=gray,linewidth=1pt, linestyle=dashed](1,2)(3,3)
\psplot[linecolor=foldline,linewidth=1pt,plotpoints=10]{1.25}{3.5}{x -2 mul 6.5 add}
\uput[0]{*0}(1.5,3.5){$\chi$}
\qdisk(1,2){2pt}
\qdisk(3,3){2pt}
\qdisk(2,2.5){2pt}
\end{pspicture}
\caption{Fold line $\chi$ places point $P$ onto line $\ell$ at $P'$.} 
\label{PbLb}
\end{figure}

Next, consider $P(p, q)$ and $\ell: x=k$ (Fig.~\ref{PbLb}). Assume that $P$ is not on $\ell$, so $p\neq k$. Fold $\chi$ places $P$ onto line $\ell$ at $P'(k, s)$, where $s$ is another parameter. A normal vector to $\chi$ is $\overrightarrow{PP'}=(k-p, s-q)$ and $\chi$ passes through the midpoint of $\overline{PP'}$ at $((s+q)/2, (k+p)/2)$. Then, $\chi$ has an equation
\begin{equation}
(k-p)x+(s-q)y=\frac{s^2-q^2}{2}+\frac{k^2-p^2}{2}.
\label{Fb}
\end{equation} 

Finally, let line $L$ be described by 
\begin{equation}
ax+by=c.
\label{L}
\end{equation}
Fold $\xi$ aligns $\chi$ and $n$, and the following two cases are possible.

\subsection{Case 1: \texorpdfstring{$\xi$ and $n$ intersect}{xi and n intersect}}

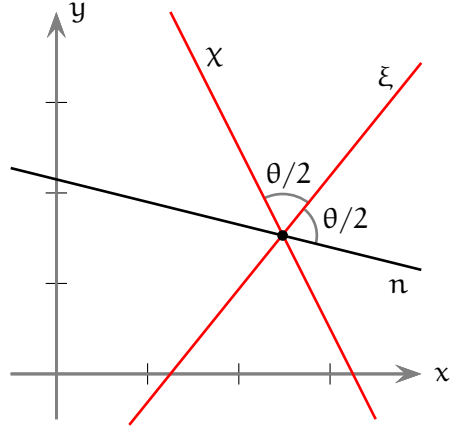
\begin{figure}[!htb]
\centering
\begin{pspicture}(-1,-1.)(5,5)
\psset{xunit=1.2cm,yunit=1.2cm}
%\showgrid
\psaxes[linecolor=gray,labels=none,linewidth=1pt]{->}(0,0)(-.5,-.5)(4,4)[$x$,0][$y$,0]
\psarc[linecolor=gray,linewidth=1pt](2.48,1.53){.45}{-15}{50}
\psarc[linecolor=gray,linewidth=1pt](2.48,1.53){.55}{50}{115}
\psplot[linecolor=foldline,linewidth=1pt,plotpoints=10]{1.25}{3.5}{x -2 mul 6.5 add}
\psplot[linecolor=foldline,linewidth=1pt,plotpoints=10]{.8}{4}{x -1.25 add 1.25 mul}
\psplot[linewidth=1pt,plotpoints=10]{-.5}{4}{x -2.48 add -.25 mul 1.53 add}
\uput[90]{*0}(3.6,3){$\xi$}
\uput[90]{*0}(3.75,0.75){$n$}
\uput[0]{*0}(1.5,3.5){$\chi$}
\uput{15pt}[15]{*0}(2.48,1.53){$\theta/2$}
\uput{17pt}[85]{*0}(2.48,1.53){$\theta/2$}
\qdisk(2.48,1.53){2pt}
\end{pspicture}
\caption{Fold line $\xi$ aligns lines $n$ and $\chi$.} 
\label{FaFb}
\end{figure}

Assume that $\xi$ and $n$ intersect (Fig.~\ref{FaFb}). Using Eqs.~(\ref{Fa}) and (\ref{L}), the intersection is found at
\begin{equation}
x=\frac{bt^2+ch}{bt+ah},\qquad y=-\frac{t(at-c)}{bt+ah},
\label{inter}
\end{equation}
provided $bt+ah\neq 0$. The intersection must be also in $\chi$; therefore, replacing in Eq.~(\ref{Fb}) we obtain
\begin{equation}
(k-p)\frac{bt^2+ch}{bt+ah}-(s-q)\frac{t(at-c)}{bt+ah}=\frac{s^2-q^2}{2}+\frac{k^2-p^2}{2}.
\label{FbInter}
\end{equation} 

Further, $\xi$ bisects the angle $\theta$ between $n$ and $\chi$. Letting $\mathbf{n}_\xi=(t, -h)$, $\mathbf{n}_\chi=(k-p, s-q)$ and $\mathbf{n}_n=(a, b)$ be the normal vectors for $\xi$, $\chi$ and $n$, respectively, then 
\begin{equation}
\cos \frac{\theta}{2}=\frac{|\mathbf{n}_\xi\cdot\mathbf{n}_\chi|}{\|\mathbf{n}_\xi\|\|\mathbf{n}_\chi\|}=\frac{|\mathbf{n}_\xi\cdot\mathbf{n}_n|}{\|\mathbf{n}_\xi\|\|\mathbf{n}_n\|},
\label{slopes}
\end{equation}
which produces
\begin{equation}
\frac{|t(k-p)-h(s-q)|}{\sqrt{(k-p)^2+(s-q)^2}}=\frac{|at-bh|}{\sqrt{a^2+b^2}}.
\label{slope2}
\end{equation}

Finally, solving Eq.~(\ref{slope2}) for $s$ and replacing in Eq.~(\ref{FbInter}), we obtain (after some algebra) a quintic equation of the form 
\begin{equation}
a^2t^5+\alpha t^4 + \beta t^3+\gamma t^2 +\delta t+\epsilon=0.
\label{quintic}
\end{equation}
 If $a=0$ then the equation becomes a quartic; this case is disregarded by letting $a=1$. Then, the equation's coefficients are
\begin{equation}
\left\{
\begin{array}{ll}
\alpha & = (-k-p+2bq-b^2k+b^2p-12bh-2c)/4,\\
\beta  &= h(q+2bp-b^2q+bc-h+2b^2h),\\
\gamma &= h^2(3p-k-6bq-b^2k-3b^2p+2bh)/2,\\
\delta & = -h^3(q+2bp-b^2q-bc),\\
\epsilon & = h^4(-k-p+2bq-b^2k+b^2p+2c)/4.
\end{array}
\right.
\label{sys}
\end{equation}
 
The next step is to solve the above system in order to obtain expressions for the parameters of $P, Q, \ell, m$ and $n$. From the first and fifth equations we obtain 
\begin{equation}
\epsilon-h^4\alpha=h^4(c+3bh),
\label{c1}
\end{equation}
and from the second and fourth,
\begin{equation}
h^2\beta+\delta=h^3(2bc+2b^2h-h).
\label{c2}
\end{equation}
Assuming  
\begin{equation}
D:=(\epsilon-h^4\alpha)^2-4h^6(h^4+h^2\beta+\delta)\ge 0
\label{D}
\end{equation}
 then Eqs.~(\ref{c1}) and (\ref{c2}) produce 
\begin{equation}
b=\frac{\epsilon -h^4\alpha \pm \sqrt{D}}{4h^5}
\label{b}
\end{equation} 
\begin{equation}
c= \frac{\epsilon - h^4\alpha \mp 3\sqrt{D}}{4h^4}
\label{c}
\end{equation}

Knowing $b$ and $c$, then a linear system of equations in $k, p$, and $q$ is obtained from Eqs.~(\ref{sys}) which has the unique solution\footnote{The existence of a unique solution may be demonstrated by following the same steps as \citet{Nishimura2015}: First, the system is written in matrix form $A\mathbf{x}=\mathbf{b}$, where $\mathbf{x}=(k, p, q)^t$. Next, reduction to row echelon form is applied to show that the ranks of both the coefficient matrix $A$ and the augmented matrix $(A | \mathbf{b})$ are equal to the number of unknowns.}
\begin{equation}
k = -\frac{17bh^3+3h^2(c+2\alpha)-\gamma}{2h^2(b^2+1)}
\label{k}
\end{equation}
\begin{equation}
p=\frac{-bh^3(b^2-3)+h^2((2\alpha-3c)b^2-c-2\alpha)+4bh\beta +(1-b^2)\gamma}{2h^2(b^2+1)^2}
\label{p}
\end{equation}
\begin{equation}
q=\frac{h^3(2b^4+4b^2+1) +bh^2(b^2c+2\alpha) +\beta h(1-b^2)-b\gamma}{h^2(b^2+1)^2}
\label{q}
\end{equation}

\subsection{Case 2: \texorpdfstring{$\xi$ and $n$ are parallel}{xi and n are parallel}}

In the case that $\xi$ and $n$ do not intersect, their respective normal vectors $(t, -h)$ and $(a, b)$ are parallel and satisfy  
\begin{equation}
bt+ah=0;
\label{d1}
\end{equation}
therefore, Eqs.~(\ref{inter}) are no longer valid. Fold $\xi$ aligns line $\chi$ and $n$, and so $\chi$ must also be parallel to $n$. Then, the normal vector of $\chi$, $(k-p, s-q)$,  is also parallel to $(a, b)$ and satisfies   
\begin{equation}
b(k-p)-a(s - q)=0.
\label{d2}
\end{equation}

Further, $\chi$ and $n$ must be equidistant from $\xi$. The distance $d$ between two parallel lines given by $\mathbf{n}\cdot\mathbf{x}=c_1$ and $\mathbf{n}\cdot\mathbf{x}=c_2$ is
\begin{equation}
d=\frac{|c_1-c_2|}{\|\mathbf{n}\|}.
\label{distance9}
\end{equation}
Using Eq.~(\ref{d1}) in Eq.~(\ref{Fa}), the equation for $\xi$ becomes
\begin{equation}
ax+by=-\frac{a^2h}{b},
\end{equation}
and therefore the distance between $\xi$ and $n$ is
\begin{equation}
d_1=\frac{|a^2h/b + c|}{\sqrt{a^2+b^2}}.
\label{distance}
\end{equation}

Similarly, using Eq.~(\ref{d2}) in Eq.~(\ref{Fb}), the equation for $\chi$ becomes
\begin{equation}
ax+by=\frac{a^2(k+p)+2abq+b^2(k-p)}{2a},
\label{tt}
\end{equation}
and therefore the distance between $\xi$ and $\chi$ is
\begin{equation}
d_2=\frac{|a^2h/b + (a^2(k+p)+2abq+b^2(k-p))/(2a)|}{\sqrt{a^2+b^2}}.
\label{distance2}
\end{equation}

Setting $d_1^2=d_2^2$ produces two conditions, 
\begin{equation}
\frac{a^2(k+p)+2abq+b^2(k-p)}{2a}=c
\label{cond1}
\end{equation}
or
\begin{equation}
4a^3h+a^2b(k+p)+2ab(bq+c)+b^3(k-p)=0
\label{cond2}
\end{equation}
Eq.~(\ref{cond1}) implies that the right hand sides of Eqs.~(\ref{L}) and (\ref{tt}) are equal and therefore $\chi=n$; thus, this condition is disregarded. Equation~(\ref{cond2}) sets a condition for obtaining parallel fold lines $\xi$ and $\chi$. However, this same equation is obtained by solving Eq.~(\ref{d1}) for $t$ and replacing into Eq.~(\ref{quintic}), and next using $a=1$ and Eqs.~(\ref{sys}) to substitute the quintic equation's coefficient. Therefore, the case of parallel $\xi$ and $n$ is contained within the solution found in the previous subsection.     

\subsection{Summary}

Given a quintic equation of the form $t^5+\alpha t^4 + \beta t^3+\gamma t^2 +\delta t+\epsilon=0,$ its solution may be found as follows: 
\begin{enumerate}
\item Find a value of $h$ that satisfies the condition in Eq.~(\ref{D})\footnote{This condition, in the case of $\alpha=0$ and $h=1$, was also obtained by \citet{Nishimura2015}, and he solved it by scaling variable $t$, as noted in the Introduction.}. If $\epsilon\neq0$, then  $D\rightarrow \epsilon^2> 0$ as $h\rightarrow 0$; therefore, the condition is always satisfied for $h$ small enough. 
If $\epsilon=0$, then $t=0$ is a solution of the quintic. Parameter $h$ defines point $Q(0, h)$ and line $m: y=-h$.

\item Use Eqs.~(\ref{b}) to (\ref{q}) obtain $b$, $c$, $k$, $p$ and $q$. These values define point $P(p, q)$, line $\ell: x=k$ and line $n: x+by=c$.  
\item  Fold along a line $\chi$ that places $P$ onto $\ell$, and simultaneously along a line $\xi$ that places $Q$ onto $m$ and aligns $n$ and $\chi$.
\item The intersection of $\xi$ with the $x$-axis is the sought value of $t$.  
\end{enumerate}
 
\section{Example}

Consider equation 
\begin{equation}
t^5+t^4-4t^3-3t^2+3t+1=0.
\label{eqt1}
\end{equation}
This equation is associated to the regular hendecagon and its solutions are \linebreak $t=2\cos(\frac{2i\pi}{11})$, with $i=1, 2,\ldots, 5$. 

\begin{figure}
\centering
\begin{pspicture}(-3,-4.)(4,2.5)
%\showgrid
\psset{xunit=.7cm,yunit=.7cm}
\psaxes[linecolor=crease,linewidth=1pt,labels=none]{->}(0,0)(-3.5,-5.5)(3.7,3)[$x$,0][$y$,0]
\uput[0](-3.5,2.7){a)}
\qdisk(0,1){2pt}
\uput[30](0,1){$Q$}
\psline[linewidth=1pt](0,-5.5)(0,3)
\uput[20](0,2.2){$L$}
\psline[linewidth=1pt](-1.5,-5.5)(-1.5,3)
\uput[180](-1.5,2.2){$\ell$}
\psline[linewidth=1pt](-3.5,-1)(4.5,-1)
\uput[90](-3,-1.1){$m$}
\psline[linecolor=crease,linewidth=1pt, linestyle=dashed](0,1)(3.37,-1)
\psplot[linewidth=1pt,linecolor=foldline,plotpoints=100]{-.7}{3}{x 1.68 mul -2.83 add}
\psplot[linewidth=1pt,linecolor=foldline,plotpoints=100]{-3.5}{4.2}{x 0.54 mul -2.83 add}
\psline[linecolor=crease,linewidth=1pt, linestyle=dashed](-2.5,-3)(-1.5, -4.84)
\qdisk(0,1){2pt}
\qdisk(3.37,-1){2pt}
\uput[-90](3.37,-1){$Q'$}
\qdisk(1.68,0){2pt}
\uput[-90](1.68,0){$t$}
\qdisk(0,-2.83){2pt}
\qdisk(-2.5,-3){2pt}
\uput[180](-2.5,-3){$P$}
\qdisk(-1.5, -4.84){2pt}
\uput[0](-1.5, -4.84){$P'$}
\qdisk(-2,-3.92){2pt}
\uput[0](1.95,2){$\xi$}
\uput[0](-3.7,-4.){$\chi$}
\end{pspicture}
\begin{pspicture}(-3,-4.)(4,2.5)
%\showgrid
\psset{xunit=.7cm,yunit=.7cm}
\psaxes[linecolor=crease,linewidth=1pt,labels=none]{->}(0,0)(-3.5,-5.5)(4.7,3)[$x$,0][$y$,0]
\uput[0](-3.5,2.7){b)}
\qdisk(0,1){2pt}
\uput[30](0,1){$Q$}
\psline[linewidth=1pt](0,-5.5)(0,3)
\uput[20](0,2.2){$L$}
\psline[linewidth=1pt](-1.5,-5.5)(-1.5,3)
\uput[180](-1.5,2.9){$\ell$}
\psline[linewidth=1pt](-3.5,-1)(4.5,-1)
\uput[90](-3,-1.1){$m$}
\psline[linecolor=crease,linewidth=1pt, linestyle=dashed](0,1)(1.67,-1)
\psplot[linewidth=1pt,linecolor=foldline,plotpoints=100]{-2}{3.5}{x .83 mul -.69 add}
\psplot[linewidth=1pt,linecolor=foldline,plotpoints=100]{-3.5}{4.2}{x -.187 mul -.7 add}
\psline[linecolor=crease,linewidth=1pt, linestyle=dashed](-2.5,-3)(-1.5, 2.35)
\qdisk(0,1){2pt}
\qdisk(1.67,-1){2pt}
\uput[-90](1.67,-1){$Q'$}
\qdisk(.83,0){2pt}
\uput[-90](.83,0){$t$}
\qdisk(0,-.69){2pt}
\qdisk(-2.5,-3){2pt}
\uput[180](-2.5,-3){$P$}
\qdisk(-1.5, 2.35){2pt}
\uput[0](-1.5, 2.35){$P'$}
\qdisk(-2,-.33){2pt}
\uput[0](2.1,2){$\xi$}
\uput[-90](3.7,-1.4){$\chi$}
\end{pspicture}

\begin{pspicture}(-3,-4.)(4,2.5)
%\showgrid
\psset{xunit=.7cm,yunit=.7cm}
\psaxes[linecolor=crease,linewidth=1pt,labels=none]{->}(0,0)(-3.5,-5.5)(4.7,3)[$x$,0][$y$,0]
\uput[0](-3.5,2.7){c)}
\qdisk(0,1){2pt}
\uput[150](0,1){$Q$}
\psline[linewidth=1pt](0,-5.5)(0,3)
\uput[20](0,2.2){$L$}
\psline[linewidth=1pt](-1.5,-5.5)(-1.5,3)
\uput[180](-1.5,2.2){$\ell$}
\psline[linewidth=1pt](-3.5,-1)(4.5,-1)
\uput[90](-3,-1.1){$m$}
\psline[linecolor=crease,linewidth=1pt, linestyle=dashed](0,1)(-.57,-1)
\psplot[linewidth=1pt,linecolor=foldline,plotpoints=100]{-3}{4}{x -.28 mul -.081 add}
\psplot[linewidth=1pt,linecolor=foldline,plotpoints=100]{-3}{1.5}{x 1.61 mul -.08 add}
\psline[linecolor=crease,linewidth=1pt, linestyle=dashed](-2.5,-3)(-1.5, -3.62)
\qdisk(0,1){2pt}
\qdisk(-.57,-1){2pt}
\uput[-90](-.57,-1){$Q'$}
\qdisk(-.28,0){2pt}
\uput[120](-.28,0){$t$}
\qdisk(0,-.081){2pt}
\qdisk(-2.5,-3){2pt}
\uput[180](-2.5,-3){$P$}
\qdisk(-1.5, -3.62){2pt}
\uput[0](-1.5, -3.62){$P'$}
\qdisk(-2,-3.31){2pt}
\uput[90](-2.5,.6){$\xi$}
\uput[0](-3.4,-4.){$\chi$}
\end{pspicture}
\begin{pspicture}(-3,-4.)(4,2.5)
%\showgrid
\psset{xunit=.7cm,yunit=.7cm}
\psaxes[linecolor=crease,linewidth=1pt,labels=none]{->}(0,0)(-3.5,-5.5)(4.7,3)[$x$,0][$y$,0]
\uput[0](-3.5,2.7){d)}
\qdisk(0,1){2pt}
\uput[30](0,1){$Q$}
\psline[linewidth=1pt](0,-5.5)(0,3)
\uput[20](0,2.2){$L$}
\psline[linewidth=1pt](-1.5,-5.5)(-1.5,3)
\uput[180](-1.5,2.2){$\ell$}
\psline[linewidth=1pt](-3.5,-1)(4.5,-1)
\uput[90](4,-1.1){$m$}
\psline[linecolor=crease,linewidth=1pt, linestyle=dashed](0,1)(-2.62,-1)
\psplot[linewidth=1pt,linecolor=foldline,plotpoints=100]{-2.5}{2}{x -1.31 mul -1.71 add}
\psplot[linewidth=1pt,linecolor=foldline,plotpoints=100]{-3.5}{3.5}{x -.27 mul -1.72 add}
\psline[linecolor=crease,linewidth=1pt, linestyle=dashed](-2.5,-3)(-1.5, .662)
\qdisk(0,1){2pt}
\qdisk(-2.62,-1){2pt}
\uput[-90](-2.62,-1){$Q'$}
\qdisk(-1.31,0){2pt}
\uput[-90](-1.31,0){$t$}
\qdisk(0,-1.71){2pt}
\qdisk(-2.5,-3){2pt}
\uput[180](-2.5,-3){$P$}
\qdisk(-1.5, .662){2pt}
\uput[60](-1.5, .662){$P'$}
\qdisk(-2,-1.17){2pt}
\uput[45](1.5,-4){$\xi$}
\uput[90](3,-2.6){$\chi$}
\end{pspicture}

\begin{pspicture}(-4,-4.)(3,2.5)
%\showgrid
\psset{xunit=.7cm,yunit=.7cm}
\psaxes[linecolor=crease,linewidth=1pt,labels=none]{->}(0,0)(-4.7,-5.5)(3.5,3)[$x$,0][$y$,0]
\uput[0](-4.7,2.7){e)}
\qdisk(0,1){2pt}
\uput[30](0,1){$Q$}
\psline[linewidth=1pt](0,-5.5)(0,3)
\uput[20](0,2.2){$L$}
\psline[linewidth=1pt](-1.5,-5.5)(-1.5,3)
\uput[180](-1.5,2.2){$\ell$}
\psline[linewidth=1pt](-4.7,-1)(3.2,-1)
\uput[90](2.7,-1.1){$m$}
\psline[linecolor=crease,linewidth=1pt, linestyle=dashed](0,1)(-3.83,-1)
\psplot[linewidth=1pt,linecolor=foldline,plotpoints=100]{-3.}{.8}{x -1.92 mul -3.68 add}
\psplot[linewidth=1pt,linecolor=foldline,plotpoints=100]{-4.5}{1.5}{x -.7 mul -3.68 add}
\psline[linecolor=crease,linewidth=1pt, linestyle=dashed](-2.5,-3)(-1.5, -1.57)
\qdisk(0,1){2pt}
\qdisk(-3.84,-1){2pt}
\uput[-90](-3.84,-1){$Q'$}
\qdisk(-1.92,0){2pt}
\uput[-90](-1.92,0){$t$}
\qdisk(0,-3.69){2pt}
\qdisk(-2.5,-3){2pt}
\uput[180](-2.5,-3){$P$}
\qdisk(-1.5, -1.57){2pt}
\uput[180](-1.5, -1.57){$P'$}
\qdisk(-2,-2.28){2pt}
\uput[180](-2.8,1.7){$\xi$}
\uput[0](1.,-4.2){$\chi$}
\end{pspicture}
\caption{Geometric solution of Eq.~(\ref{eqt1}), where $t=2\cos(2i\pi/11)$ and a) $i=1$, b) $i=2$, c) $i=3$, d) $i=4$, e) $i=5$.} 
\label{twofold}
\end{figure}
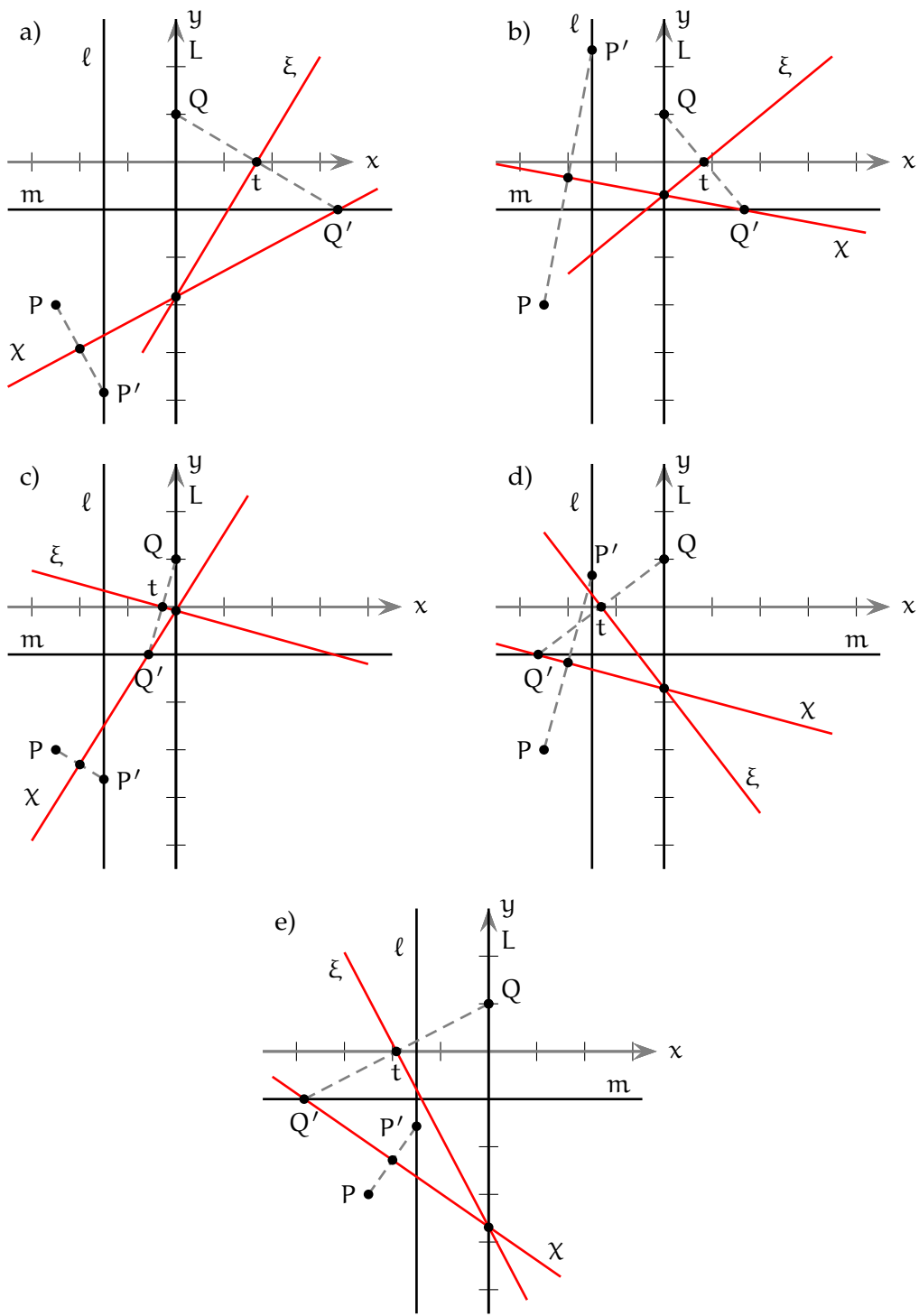

A value of $h=1$ in Eq.~(\ref{D}) produces $D=0$, and Eqs.~(\ref{b})  to (\ref{q}) result in $b=c=0$, $k=-\frac{3}{2}$, $p=-\frac{5}{2}$, and $q=-3$. The five possible geometric solutions of Eq.~(\ref{eqt1}) are illustrated in Fig.~\ref{twofold}.

As a comparison, let us note that direct application of \citeauthor{Nishimura2015}'s \citeyearpar{Nishimura2015} analysis to the same equation requires first a reduction to the depressed form $t^5-\frac{22}{5}t^3-\frac{11}{25}t^2+\frac{462}{125}t+\frac{979}{3125}=0$. This polynomial does not satisfy $\epsilon^2-4(\beta + \delta +1)\ge 0$, and so a scaling factor must be applied. For example, $t=t'/5$ produces $t'^5-110t^3-55t^2+2310t+979=0$ which satisfies the previous condition. For this polynomial, \citeauthor{Nishimura2015}'s solution is $Q(0, 1)$, $m: y=-1$; and the remaining parameters may be obtained from  Eqs.~(\ref{D})  to (\ref{q}) with $h=1$, which produce  $L:= (\frac{979\pm 3\sqrt{D}}{1897073})x+(\frac{951838\pm 979\sqrt{D}}{1897073})y+1=0$, $P(\frac{-18847612471337914\pm 19341103112746\sqrt{D}}{75005524688329}, \frac{10933226315913024\mp 11142037150475\sqrt{D}}{75005524688329})$, $\ell: x=\frac{-2758198234\pm2830266\sqrt{D}}{8660573}$, where $D=949637$. Naturally, this set of points and lines is much more difficult to produce by paper folding than the result given above.

%\bibliographystyle{apalike}
%\bibliography{../../origami}

\end{document}